\documentstyle[11pt]{article}
\begin{document}
\newcommand{\p}{\parallel }
\makeatletter \makeatother
\newtheorem{th}{Theorem}[section]
\newtheorem{lem}{Lemma}[section]
\newtheorem{de}{Definition}[section]
\newtheorem{rem}{Remark}[section]
\newtheorem{cor}{Corollary}[section]
\renewcommand{\theequation}{\thesection.\arabic {equation}}

\title{{\bf A note on modular forms and generalized anomaly cancellation formulas}}

\author{Kefeng Liu and Yong Wang\\
 }

\date{}
\maketitle

\begin{abstract}~~  By studying modular invariance
properties of some characteristic forms, we prove some new
anomaly cancellation formulas which generalize the Han-Zhang and Han-Liu-Zhang anomaly cancellation formulas \\

 \noindent{\bf Keywords:}\quad
 Modular invariance; Anomaly cancellation formulas; \\
\end{abstract}

\section{Introduction}
  \quad In 1983, the physicists Alvarez-Gaum\'{e} and Witten [AW]
  discovered the "miraculous cancellation" formula for gravitational
  anomaly which reveals a beautiful relation between the top
  components of the Hirzebruch $\widehat{L}$-form and
  $\widehat{A}$-form of a $12$-dimensional smooth Riemannian
  manifold. Kefeng Liu [L] established higher dimensional "miraculous cancellation"
  formulas for $(8k+4)$-dimensional Riemannian manifolds by
  developing modular invariance properties of characteristic forms.
  These formulas could be used to deduce some divisibility results. In
  [HZ1], [HZ2], for each $(8k+4)$-dimensional smooth Riemannian
  manifold, a more general cancellation formula that involves a
  complex line bundle was established. This formula was applied to
  ${\rm spin}^c$ manifolds, then an analytic Ochanine congruence
  formula was derived. In [CH1], Qingtao Chen and Fei Han obtained more twisted cancellation formulas for $8k$ and $8k+4$ dimensional
  manifolds and they also applied their cancellation formulas to study divisibilities on spin manifolds and
  congruences on ${\rm spin}^c$ manifolds. Recently, Han, Liu and
  Zhang generalized the anomaly cancellation formulas to cases that
  an auxiliary bundle as well as a complex line bundle are involved
  with non conditions on the first Pontryagin forms being assumed in [HLZ].
 On the other hand, motivated by the Chern-Simons theory, in
  [CH2], Qingtao Chen and Fei Han computed the transgressed forms of some modularly invariant characteristic forms,
which are related to the elliptic genera. They studied the
modularity properties of these secondary characteristic forms and
relations among them. In [W], Wang generalized the Chen-Han
cancellation
  formulas to the case that a complex line bundle is involved.
  One naturally asks if there exist more cancellation
  formulas similar to the Han-Liu-Zhang cancellation formulas. In
  this note, we derive more these type cancellation formulas.

\indent This paper is organized as follows: In Section 2, we review
some knowledge on characteristic forms and modular forms that we are
going to use. In Section 3, we generalize the Han-Liu-Zhang
cancellation formulas to the (a,b) type cancellation formulas. In
Section 4, we prove some general cancellation formulas involving two
  complex line bundles. In Section 5, we study modular invariance
properties of some characteristic forms on odd dimensional manifolds.\\

\section{characteristic forms and modular forms}
 \quad The purpose of this section is to review the necessary knowledge on
characteristic forms and modular forms that we are going to use.\\

 \noindent {\bf  2.1 characteristic forms. }Let $M$ be a Riemannian manifold.
 Let $\nabla^{ TM}$ be the associated Levi-Civita connection on $TM$
 and $R^{TM}=(\nabla^{TM})^2$ be the curvature of $\nabla^{ TM}$.
 Let $\widehat{A}(TM,\nabla^{ TM})$ and $\widehat{L}(TM,\nabla^{ TM})$
 be the Hirzebruch characteristic forms defined respectively by (cf.
 [Z])
 $$\widehat{A}(TM,\nabla^{ TM})={\rm
 det}^{\frac{1}{2}}\left(\frac{\frac{\sqrt{-1}}{4\pi}R^{TM}}{{\rm
 sinh}(\frac{\sqrt{-1}}{4\pi}R^{TM})}\right),$$
 $$\widehat{L}(TM,\nabla^{ TM})={\rm
 det}^{\frac{1}{2}}\left(\frac{\frac{\sqrt{-1}}{2\pi}R^{TM}}{{\rm
 tanh}(\frac{\sqrt{-1}}{4\pi}R^{TM})}\right).\eqno(2.1)$$
   Let $E$, $F$ be two Hermitian vector bundles over $M$ carrying
   Hermitian connection $\nabla^E,\nabla^F$ respectively. Let
   $R^E=(\nabla^E)^2$ (resp. $R^F=(\nabla^F)^2$) be the curvature of
   $\nabla^E$ (resp. $\nabla^F$). If we set the formal difference
   $G=E-F$, then $G$ carries an induced Hermitian connection
   $\nabla^G$ in an obvious sense. We define the associated Chern
   character form as
   $${\rm ch}(G,\nabla^G)={\rm tr}\left[{\rm
   exp}(\frac{\sqrt{-1}}{2\pi}R^E)\right]-{\rm tr}\left[{\rm
   exp}(\frac{\sqrt{-1}}{2\pi}R^F)\right].\eqno(2.2)$$
   For any complex number $t$, let
   $$\wedge_t(E)={\bf C}|_M+tE+t^2\wedge^2(E)+\cdots,~S_t(E)={\bf
   C}|_M+tE+t^2S^2(E)+\cdots$$
   denote respectively the total exterior and symmetric powers of
   $E$, which live in $K(M)[[t]].$ The following relations between
   these operations hold,
   $$S_t(E)=\frac{1}{\wedge_{-t}(E)},~\wedge_t(E-F)=\frac{\wedge_t(E)}{\wedge_t(F)}.\eqno(2.3)$$
   Moreover, if $\{\omega_i\},\{\omega_j'\}$ are formal Chern roots
   for Hermitian vector bundles $E,F$ respectively, then
   $${\rm ch}(\wedge_t(E))=\prod_i(1+e^{\omega_i}t).\eqno(2.4)$$
   Then we have the following formulas for Chern character forms,
   $${\rm ch}(S_t(E))=\frac{1}{\prod_i(1-e^{\omega_i}t)},~
{\rm
ch}(\wedge_t(E-F))=\frac{\prod_i(1+e^{\omega_i}t)}{\prod_j(1+e^{\omega_j'}t)}.\eqno(2.5)$$
\indent If $W$ is a real Euclidean vector bundle over $M$ carrying a
Euclidean connection $\nabla^W$, then its complexification $W_{\bf
C}=W\otimes {\bf C}$ is a complex vector bundle over $M$ carrying a
canonical induced Hermitian metric from that of $W$, as well as a
Hermitian connection $\nabla^{W_{\bf C}}$ induced from $\nabla^W$.
If $E$ is a vector bundle (complex or real) over $M$, set
$\widetilde{E}=E-{\rm dim}E$ in $K(M)$ or $KO(M)$.\\

\noindent{\bf 2.2 Some properties about the Jacobi theta functions
and modular forms}\\
   \indent We first recall the four Jacobi theta functions are
   defined as follows( cf. [C]):
   $$\theta(v,\tau)=2q^{\frac{1}{8}}{\rm sin}(\pi
   v)\prod_{j=1}^{\infty}[(1-q^j)(1-e^{2\pi\sqrt{-1}v}q^j)(1-e^{-2\pi\sqrt{-1}v}q^j)],\eqno(2.6)$$
$$\theta_1(v,\tau)=2q^{\frac{1}{8}}{\rm cos}(\pi
   v)\prod_{j=1}^{\infty}[(1-q^j)(1+e^{2\pi\sqrt{-1}v}q^j)(1+e^{-2\pi\sqrt{-1}v}q^j)],\eqno(2.7)$$
$$\theta_2(v,\tau)=\prod_{j=1}^{\infty}[(1-q^j)(1-e^{2\pi\sqrt{-1}v}q^{j-\frac{1}{2}})
(1-e^{-2\pi\sqrt{-1}v}q^{j-\frac{1}{2}})],\eqno(2.8)$$
$$\theta_3(v,\tau)=\prod_{j=1}^{\infty}[(1-q^j)(1+e^{2\pi\sqrt{-1}v}q^{j-\frac{1}{2}})
(1+e^{-2\pi\sqrt{-1}v}q^{j-\frac{1}{2}})],\eqno(2.9)$$ \noindent
where $q=e^{2\pi\sqrt{-1}\tau}$ with $\tau\in\textbf{H}$, the upper
half complex plane. Let
$$\theta'(0,\tau)=\frac{\partial\theta(v,\tau)}{\partial
v}|_{v=0}.\eqno(2.10)$$ \noindent Then the following Jacobi identity
(cf. [Ch]) holds,
$$\theta'(0,\tau)=\pi\theta_1(0,\tau)\theta_2(0,\tau)\theta_3(0,\tau).\eqno(2.11)$$
\noindent Denote $SL_2({\bf Z})=\left\{\left(\begin{array}{cc}
\ a & b  \\
 c  & d
\end{array}\right)\mid a,b,c,d \in {\bf Z},~ad-bc=1\right\}$ the
modular group. Let $S=\left(\begin{array}{cc}
\ 0 & -1  \\
 1  & 0
\end{array}\right),~T=\left(\begin{array}{cc}
\ 1 &  1 \\
 0  & 1
\end{array}\right)$ be the two generators of $SL_2(\bf{Z})$. They
act on $\textbf{H}$ by $S\tau=-\frac{1}{\tau},~T\tau=\tau+1$. One
has the following transformation laws of theta functions under the
actions of $S$ and $T$ (cf. [C]):
$$\theta(v,\tau+1)=e^{\frac{\pi\sqrt{-1}}{4}}\theta(v,\tau),~~\theta(v,-\frac{1}{\tau})
=\frac{1}{\sqrt{-1}}\left(\frac{\tau}{\sqrt{-1}}\right)^{\frac{1}{2}}e^{\pi\sqrt{-1}\tau
v^2}\theta(\tau v,\tau);\eqno(2.12)$$
$$\theta_1(v,\tau+1)=e^{\frac{\pi\sqrt{-1}}{4}}\theta_1(v,\tau),~~\theta_1(v,-\frac{1}{\tau})
=\left(\frac{\tau}{\sqrt{-1}}\right)^{\frac{1}{2}}e^{\pi\sqrt{-1}\tau
v^2}\theta_2(\tau v,\tau);\eqno(2.13)$$
$$\theta_2(v,\tau+1)=\theta_3(v,\tau),~~\theta_2(v,-\frac{1}{\tau})
=\left(\frac{\tau}{\sqrt{-1}}\right)^{\frac{1}{2}}e^{\pi\sqrt{-1}\tau
v^2}\theta_1(\tau v,\tau);\eqno(2.14)$$
$$\theta_3(v,\tau+1)=\theta_2(v,\tau),~~\theta_3(v,-\frac{1}{\tau})
=\left(\frac{\tau}{\sqrt{-1}}\right)^{\frac{1}{2}}e^{\pi\sqrt{-1}\tau
v^2}\theta_3(\tau v,\tau),\eqno(2.15)$$
 $$\theta'(v,\tau+1)=e^{\frac{\pi\sqrt{-1}}{4}}\theta'(v,\tau),~~
 \theta'(0,-\frac{1}{\tau})=\frac{1}{\sqrt{-1}}\left(\frac{\tau}{\sqrt{-1}}\right)^{\frac{1}{2}}
\tau\theta'(0,\tau).\eqno(2.16)$$
 \noindent {\bf Definition 2.1} A modular form over $\Gamma$, a
 subgroup of $SL_2({\bf Z})$, is a holomorphic function $f(\tau)$ on
 $\textbf{H}$ such that
 $$f(g\tau):=f\left(\frac{a\tau+b}{c\tau+d}\right)=\chi(g)(c\tau+d)^kf(\tau),
 ~~\forall g=\left(\begin{array}{cc}
\ a & b  \\
 c & d
\end{array}\right)\in\Gamma,\eqno(2.17)$$
\noindent where $\chi:\Gamma\rightarrow {\bf C}^{\star}$ is a
character of $\Gamma$. $k$ is called the weight of $f$.\\
Let $$\Gamma_0(2)=\left\{\left(\begin{array}{cc}
\ a & b  \\
 c  & d
\end{array}\right)\in SL_2({\bf Z})\mid c\equiv 0~({\rm
mod}~2)\right\},$$
$$\Gamma^0(2)=\left\{\left(\begin{array}{cc}
\ a & b  \\
 c  & d
\end{array}\right)\in SL_2({\bf Z})\mid b\equiv 0~({\rm
mod}~2)\right\},$$  be the two modular subgroups of $SL_2({\bf Z})$.
It is known that the generators of $\Gamma_0(2)$ are
$T,~ST^2ST$, the generators of $\Gamma^0(2)$ are $STS,~T^2STS$ (cf.[C]).\\
\indent Let $E_2(\tau)$ be Eisenstein series which is a quasimodular
form over $SL(2,{\bf Z})$, satisfying:
$$E_2\left(\frac{a\tau +b}{c\tau
+d}\right)=(c\tau+d)^2E_2(\tau)-\frac{6\sqrt{-1}c(c\tau+d)}{\pi}.\eqno(2.18)$$
In particular, we have
$$E_2(\tau+1)=E_2(\tau),\eqno(2.19)$$
$$E_2(-\frac{1}{\tau})=\tau^2E_2(\tau)-\frac{6\sqrt{-1}\tau}{\pi}.\eqno(2.20)$$
\indent If $\Gamma$ is a modular subgroup, let ${\mathcal{M}}_{{\bf
R}}(\Gamma)$ denote the ring of modular forms over $\Gamma$ with
real Fourier coefficients. Writing $\theta_j=\theta_j(0,\tau),~1\leq
j\leq 3,$ we introduce four explicit modular forms (cf. [L]),
$$\delta_1(\tau)=\frac{1}{8}(\theta_2^4+\theta_3^4),~~\varepsilon_1(\tau)=\frac{1}{16}\theta_2^4\theta_3^4,$$
$$\delta_2(\tau)=-\frac{1}{8}(\theta_1^4+\theta_3^4),~~\varepsilon_2(\tau)=\frac{1}{16}\theta_1^4\theta_3^4.$$
\noindent They have the following Fourier expansions in
$q^{\frac{1}{2}}$:
$$\delta_1(\tau)=\frac{1}{4}+6q+\cdots,~~\varepsilon_1(\tau)=\frac{1}{16}-q+\cdots,$$
$$\delta_2(\tau)=-\frac{1}{8}-3q^{\frac{1}{2}}+\cdots,~~\varepsilon_2(\tau)=q^{\frac{1}{2}}+\cdots,$$
\noindent where the $"\cdots"$ terms are the higher degree terms,
all of which have integral coefficients. They also satisfy the
transformation laws,
$$\delta_2(-\frac{1}{\tau})=\tau^2\delta_1(\tau),~~~~~~\varepsilon_2(-\frac{1}{\tau})
=\tau^4\varepsilon_1(\tau),\eqno(2.21)$$

\noindent {\bf Lemma 2.2} ([L]) {\it $\delta_1(\tau)$ (resp.
$\varepsilon_1(\tau)$) is a modular form of weight $2$ (resp. $4$)
over $\Gamma_0(2)$, $\delta_2(\tau)$ (resp. $\varepsilon_2(\tau)$)
is a modular form of weight $2$ (resp. $4$) over $\Gamma^0(2)$ and
moreover ${\mathcal{M}}_{{\bf R}}(\Gamma^0(2))={\bf
R}[\delta_2(\tau),\varepsilon_2(\tau)]$.}

\section { A generalization of the Han-Liu-Zhang cancellation formulas }

\quad Let $M$ be a $4k$ dimensional Riemannian manifold and $V$ be a
rank $2l$ real vector bundle on $M$. Let $a,~b$ be two integers.
  Set
   $$\Theta_1(T_{C}M,V_C,a,b)=
   \bigotimes _{n=1}^{\infty}S_{q^n}(\widetilde{T_CM})\otimes
\bigotimes _{m=1}^{\infty}\wedge_{q^m}(\widetilde{V_C})^a$$
$$~~~~~~~~\otimes (\bigotimes _{r=1}^{\infty}\wedge
_{q^{r-\frac{1}{2}}}(\widetilde{V_C}))^b\otimes(\bigotimes
_{s=1}^{\infty}\wedge _{-q^{s-\frac{1}{2}}}(\widetilde{V_C}))^b,$$
$$\Theta_2(T_{C}M,V_C,a,b)=\bigotimes _{n=1}^{\infty}S_{q^n}(\widetilde{T_CM})\otimes
(\bigotimes _{m=1}^{\infty}\wedge_{q^{m}}(\widetilde{V_C}))^b$$
$$~~~~~~~~\otimes (\bigotimes _{r=1}^{\infty}\wedge
_{q^{r-\frac{1}{2}}}(\widetilde{V_C}))^b\otimes(\bigotimes
_{s=1}^{\infty}\wedge
_{-q^{s-\frac{1}{2}}}(\widetilde{V_C}))^a,\eqno(3.1)$$ Clearly,
$\Theta_1(T_{C}M,V_C,a,b)$ and $\Theta_2(T_{C}M,V_C,a,b)$ admit
formal Fourier expansion in $q^{\frac{1}{2}}$ as
$$\Theta_1(T_{C}M,V_C,a,b)=A_0(T_{C}M,V_C,a,b)+A_1(T_{C}M,V_C,a,b)q
^{\frac{1}{2}}+\cdots,$$
$$\Theta_2(T_{C}M,V_C,a,b)=B_0(T_{C}M,V_C,a,b)+B_1(T_{C}M,V_C,a,b)q
^{\frac{1}{2}}+\cdots,\eqno(3.2)$$ where the $A_j$ and $B_j$ are
elements in the semi-group formally generated by Hermitian vector
bundles over $M$. Moreover, they carry canonically induced Hermitian
connections. Let $\{\pm 2\pi iy_\nu\}$ be the formal Chern roots of
$V_C$. If $V$ is spin and $\triangle(V)$ is the spinor bundle of
$V$, one know that the Chern character of $\triangle(V)$ is given by
$${\rm ch}(\triangle(V))=\prod_{\nu=1}^l(e^{\pi iy_\nu}+e^{-\pi
iy_\nu}).$$
 In the following, we do not assume that $V$ is spin, but still
 formally use ${\rm ch}((\triangle(V))^a)$ for the short-hand notion
 of $(\prod_{\nu=1}^l(e^{\pi iy_\nu}+e^{-\pi
iy_\nu}))^a$ which is a well-defined cohomology class on $M$. Let
$p_1$ denote the first Pontryagin class.
  If $\omega$ is a differential form over $M$, we denote
$\omega^{(4k)}$ its top degree component. Define virtual complex
vector bundle $b_r(T_{\bf C}M, V_{\bf C},a,b)$ on $M$, $0\leq r\leq
[\frac{k}{2}],$ via the equality
$$\Theta_2(T_{C}M,V_C,a,b)\equiv \sum
_{r=0}^{[\frac{k}{2}]}b_r(8\delta_2)^{k-2r}\varepsilon_2^r~~~~{\rm
mod}~q^{\frac{[\frac{k}{2}]+1}{2}}\dot
K(M)[[q^{\frac{1}{2}}]].\eqno(3.3)$$ Then
$$b_0=(-1)^k{\bf C},~~b_1=-24(-1)^kk-a\widetilde{V_{\bf
C}}.\eqno(3.4)$$
 \indent Define degree $4k-4$ differential forms $\beta_r(T_{\bf C}M, V_{\bf C},a,b)$ on $M$, $0\leq r\leq
[\frac{k}{2}],$ via the equality
$$\left\{\frac{e^{\frac{1}{24}E_2(\tau)[p_1(TM)-(a+2b)p_1(V)]}-1}{p_1(TM)-(a+2b)p_1(V)}\widehat{A}(TM)
{\rm ch}((\triangle(V))^b){\rm
ch}(\Theta_2(T_{C}M,V_C,a,b))\right\}^{(4k-4)}$$ $$\equiv \sum
_{r=0}^{[\frac{k}{2}]}\beta_r(8\delta_2)^{k-2r}\varepsilon_2^r~~~~{\rm
mod}~q^{\frac{[\frac{k}{2}]+1}{2}}\dot
\Omega^{(4k-4)}(M)[[q^{\frac{1}{2}}]].\eqno(3.5)$$ It is easy to
calculate that
$$\beta_0=(-1)^k\left\{\frac{e^{\frac{1}{24}[p_1(TM)-(a+2b)p_1(V)]}-1}{p_1(TM)-(a+2b)p_1(V)}\widehat{A}(TM)
{\rm ch}((\triangle(V))^b)\right\}^{(4k-4)},\eqno(3.6)$$
$$\beta_1=(-1)^k\left\{\frac{e^{\frac{1}{24}[p_1(TM)-(a+2b)p_1(V)]}-1}{p_1(TM)-(a+2b)p_1(V)}\widehat{A}(TM)
{\rm ch}((\triangle(V))^b){\rm ch}(-a\widetilde{V_{\bf
C}}-24k)\right\}^{(4k-4)}.\eqno(3.7)$$
Our main results in this section include the following theorem.\\

\noindent {\bf Theorem 3.1}
$$\left\{{\widehat{A}(TM,\nabla^{TM})}{\rm
ch}((\triangle(V))^a)\right\}^{(4k)}$$
$$-\sum_{r=0}^{[\frac{k}{2}]}2^{(a-b)l+k-6r}\left\{\widehat{A}(TM,\nabla^{TM})
{\rm ch}((\triangle(V))^b){\rm ch}(b_r)\right\}^{(4k)}$$
$$=[p_1(TM)-(a+2b)p_1(V)]{\cal
B}(\nabla^{TM},\nabla^V,a,b),\eqno(3.8)$$
 {\it where}
$${\cal
B}(\nabla^{TM},\nabla^V,a,b)=\sum_{r=0}^{[\frac{k}{2}]}2^{(a-b)l+k-6r}\beta_r$$
$$-\left\{\frac{e^{\frac{1}{24}[p_1(TM)-(a+2b)p_1(V)]}-1}{p_1(TM)-(a+2b)p_1(V)}\widehat{A}(TM)
{\rm ch}((\triangle(V))^a)\right\}^{(4k-4)}.\eqno(3.9)$$

\noindent {\bf Proof.} Let $\{\pm2\pi\sqrt{-1}x_j| ~1\leq j\leq
2k\}$
 be the Chern roots of $T_CM$ .
 Set
$$Q_1(\tau)=\left\{e^{\frac{1}{24}E_2(\tau)[p_1(TM)-(a+2b)p_1(V)]}{\widehat{A}(TM,\nabla^{TM})}{\rm ch}((\triangle(V))^a)
 {\rm
ch}(\Theta_1(T_CM,V_C,a,b))\right\}^{(4k)},\eqno(3.10)$$
$$Q_2(\tau)=\left\{{\widehat{A}(TM,\nabla^{TM})}{\rm ch}((\triangle(V))^b)
 {\rm
ch}(\Theta_2(T_CM,V_C,a,b))\right\}^{(4k)},\eqno(3.11)$$
$$\overline{Q_2(\tau)}=\left\{\frac{e^{\frac{1}{24}E_2(\tau)[p_1(TM)-(a+2b)p_1(V)]}-1}{p_1(TM)-(a+2b)p_1(V)}
\right.$$ $$ \left.\dot{\widehat{A}(TM,\nabla^{TM})}{\rm
ch}((\triangle(V))^b)
 {\rm
ch}(\Theta_2(T_CM,V_C,a,b))\right\}^{(4k-4)},\eqno(3.12)$$
 Direct
computations show that
$$Q_1(\tau)=2^{al}\left\{e^{\frac{1}{24}E_2(\tau)[p_1(TM)-(a+2b)p_1(V)]}
\left(\prod_{j=1}^{2k}\frac{x_j\theta'(0,\tau)}{\theta(x_j,\tau)}\right)\right.$$
$$\left. \cdot \left(\prod_{\nu=1}^{l}\frac{\theta^a_1(y_\nu,\tau)}{\theta^a_1(0,\tau)}
\frac{\theta^b_2(y_\nu,\tau)}{\theta^b_2(0,\tau)}
\frac{\theta^b_3(y_\nu,\tau)}{\theta^b_3(0,\tau)}\right)\right\}^{(4k)}.\eqno(3.13)$$
Similarly, we have
$$Q_2(\tau)+[p_1(TM)-(a+2b)p_1(V)]\overline{Q_2(\tau)}$$
$$=2^{bl}\left\{e^{\frac{1}{24}E_2(\tau)[p_1(TM)-(a+2b)p_1(V)]}
\left(\prod_{j=1}^{2k}\frac{x_j\theta'(0,\tau)}{\theta(x_j,\tau)}\right)\right.$$
$$\left. \cdot\left(\prod_{\nu=1}^{l}\frac{\theta^a_2(y_\nu,\tau)}{\theta^a_2(0,\tau)}
\frac{\theta^b_1(y_\nu,\tau)}{\theta^b_1(0,\tau)}
\frac{\theta^b_3(y_\nu,\tau)}{\theta^b_3(0,\tau)}\right)\right\}^{(4k)}.\eqno(3.14)$$
  By
(2.12)-(2.16) and (2.19)-(2.20), then $Q_1(\tau)$ is a modular form
of weight $2k$ over $\Gamma_0(2)$, while
$Q_2(\tau)+[p_1(TM)-(a+2b)p_1(V)]\overline{Q_2(\tau)}$ is a modular
form of weight $2k$ over $\Gamma^0(2)$ . Moreover, the following
identity holds,
$$Q_1(-\frac{1}{\tau})=2^{(a-b)l}\tau^{2k}(Q_2(\tau)+[p_1(TM)-(a+2b)p_1(V)]\overline{Q_2(\tau)}).\eqno(3.15)$$
\indent Observe that at any point $x\in M$, up to the volume form
determined by the metric on $T_xM$, both $Q_i(\tau),$ and
$Q_2(\tau)+[p_1(TM)-(a+2b)p_1(V)]\overline{Q_2(\tau)}$ can be view
as a power series of $q^{\frac{1}{2}}$ with real Fourier
coefficients. By Lemma 2.2, we have
$$Q_2(\tau)+[p_1(TM)-(a+2b)p_1(V)]\overline{Q_2(\tau)}=h_0(8\delta_2)^{k}+h_1(8\delta_2)^{k-2}\varepsilon_2+\cdots+h_{[\frac{k}{2}]}(8\delta_2)
^{k-2[\frac{k}{2}]}\varepsilon^{[\frac{k}{2}]}_2,\eqno(3.16)$$ where
each $h_r,~ 0\leq r\leq [\frac{k}{2}],$ is a real multiple of the
volume form at $x$. By the definitions of $b_r$ and $\beta_r$, we
have
$$h_r=\left\{\widehat{A}(TM,\nabla^{TM})
{\rm ch}((\triangle(V))^b){\rm
ch}(b_r)\right\}^{(4k)}+[p_1(TM)-(a+2b)p_1(V)]\beta_r.\eqno(3.17)$$
By (2.21) (3.15) and (3.16), we get
$$Q_1(\tau)=2^{(a-b)l}\left[h_0(8\delta_1)^{k}+h_1(8\delta_1)^{k-2}\varepsilon_1+\cdots+h_{[\frac{k}{2}]}(8\delta_1)
^{k-2[\frac{k}{2}]}\varepsilon^{[\frac{k}{2}]}_1\right].\eqno(3.18)$$
By comparing the constant term in (3.18), we get
$$ \left\{e^{\frac{1}{24}E_2(\tau)[p_1(TM)-(a+2b)p_1(V)]}\widehat{A}(TM,\nabla^{TM})
{\rm ch}((\triangle(V))^a)\right\}^{(4k)}$$ $$=
\sum_{r=0}^{[\frac{k}{2}]}2^{(a-b)l+k-6r}\left(\left\{\widehat{A}(TM,\nabla^{TM})
{\rm ch}((\triangle(V))^b){\rm ch}(b_r)\right\}^{(4k)}
+[p_1(TM)-(a+2b)p_1(V)]\beta_r\right),\eqno(3.19)$$ So we get
Theorem 3.1.
$\Box$ \\

\noindent {\bf Remark.} By Theorem 3.1, if $M$ is a $4k$ dimensional
spin manifolds and $V$ is a $2l$ dimensional spin vector bundle over
$M$, then we have
$$\int_M[p_1(TM)-(a+2b)p_1(V)]{\cal
B}(\nabla^{TM},\nabla^V,a,b),\eqno(3.20)$$ is an integer. When $a=1,
~b=0$, we get the Han-Liu-Zhang cancellation formula.

By the direct computations, we have \\

\noindent {\bf Corollary 3.2} {\it When ${\rm dim}M=4$, the
following identity holds}
$$\left\{\widehat{A}(TM,\nabla^{TM}) {\rm
ch}((\triangle(V))^a)\right\}^{(4)}+2^{(a-b)l+1}\left\{\widehat{A}(TM,\nabla^{TM})
{\rm ch}((\triangle(V))^b)\right\}^{(4)}$$
$$=-2^{al-3}[p_1(TM)-(a+2b)p_1(V)].\eqno(3.20)$$

\noindent {\bf Corollary 3.3} {\it When ${\rm dim}M=8$, the
following identity holds}
$$\left\{\widehat{A}(TM,\nabla^{TM}) {\rm
ch}((\triangle(V))^a)\right\}^{(8)}$$
$$-\left\{-2^{(a-b)l-4}a
\widehat{A}(TM,\nabla^{TM}) {\rm ch}((\triangle(V))^b){\rm
ch}(\widetilde{V_{\bf C}})+2^{(a-b)l}\widehat{A}(TM,\nabla^{TM})
{\rm ch}((\triangle(V))^b)\right\}^{(8)} $$ $$
=[p_1(TM)-(a+2b)p_1(V)]
\left\{\frac{e^{\frac{1}{24}[p_1(TM)-(a+2b)p_1(V)]}-1}{p_1(TM)-(a+2b)p_1(V)}\right.$$
$$\left.
[2^{(a-b)l}\widehat{A}(TM,\nabla^{TM}) {\rm ch}((\triangle(V))^b)-a
\widehat{A}(TM,\nabla^{TM}) {\rm ch}((\triangle(V))^b){\rm
ch}(\widetilde{V_{\bf C}})\right.$$
$$\left.-\widehat{A}(TM,\nabla^{TM})
{\rm ch}((\triangle(V))^a)]\right\}^{(8)}.\eqno(3.21)$$

   \quad Let $\xi$ be a rank two real oriented Euclidean vector
   bundle over $M$ carrying with an Euclidean connection
   $\nabla^{\xi}$. Set
   $$\Theta_1(T_{C}M,V_C,\xi_C,a,b)=
   \bigotimes _{n=1}^{\infty}S_{q^n}(\widetilde{T_CM})\otimes
\bigotimes
_{m=1}^{\infty}\wedge_{q^m}(a\widetilde{V_C}-2\widetilde{\xi_C})$$
$$~~~~~~~~\otimes \bigotimes _{r=1}^{\infty}\wedge
_{q^{r-\frac{1}{2}}}(b\widetilde{V_C}+\widetilde{\xi_C})\otimes\bigotimes
_{s=1}^{\infty}\wedge
_{-q^{s-\frac{1}{2}}}(b\widetilde{V_C}+\widetilde{\xi_C}),$$
$$\Theta_2(T_{C}M,V_C,\xi_C,a,b)=\bigotimes _{n=1}^{\infty}S_{q^n}(\widetilde{T_CM})\otimes
\bigotimes
_{m=1}^{\infty}\wedge_{q^{m}}(b\widetilde{V_C}+\widetilde{\xi_C})$$
$$~~~~~~~~\otimes \bigotimes _{r=1}^{\infty}\wedge
_{q^{r-\frac{1}{2}}}(b\widetilde{V_C}+\widetilde{\xi_C})\otimes\bigotimes
_{s=1}^{\infty}\wedge
_{-q^{s-\frac{1}{2}}}(a\widetilde{V_C}-2\widetilde{\xi_C}),\eqno(3.22)$$
Clearly, $\Theta_1(T_{C}M,V_C,\xi_C,a,b)$ and
$\Theta_2(T_{C}M,V_C,\xi_C,a,b)$ admit formal Fourier expansion in
$q^{\frac{1}{2}}$ as
$$\Theta_1(T_{C}M,V_C,\xi_C,a,b)=A_0(T_{C}M,V_C,\xi_C,a,b)+A_1(T_{C}M,V_C,\xi_C,a,b)q
^{\frac{1}{2}}+\cdots,$$
$$\Theta_2(T_{C}M,V_C,\xi_C,a,b)=B_0(T_{C}M,V_C,\xi_C,a,b)+B_1(T_{C}M,V_C,\xi_C,a,b)q
^{\frac{1}{2}}+\cdots,\eqno(3.23)$$ Let $c=2\pi\sqrt{-1} u$ be the
Euler form of $\xi$. Set
$$\widetilde{Q}_1(\tau)=\left\{
\frac{e^{\frac{1}{24}E_2(\tau)[p_1(TM)-(a+2b)p_1(V)]}{\widehat{A}(TM,\nabla^{TM})}{\rm
ch}((\triangle(V))^a)}{{\rm cosh}^2(\frac{c}{2})}
 {\rm
ch}(\Theta_1(T_CM,V_C,\xi_C,a,b))\right\}^{(4k)},\eqno(3.24)$$
$$\widetilde{Q}_2(\tau)=\left\{{\widehat{A}(TM,\nabla^{TM})}{\rm cosh}(\frac{c}{2}){\rm ch}((\triangle(V))^b)
 {\rm
ch}(\Theta_2(T_CM,V_C,\xi_C,a,b))\right\}^{(4k)},\eqno(3.25)$$
$$\widetilde{Q}_3(\tau)=\left\{\frac{e^{\frac{1}{24}E_2(\tau)[p_1(TM)-(a+2b)p_1(V)]}-1}{p_1(TM)-(a+2b)p_1(V)}
\right.$$
$$\left.\cdot{\widehat{A}(TM,\nabla^{TM})}{\rm cosh}(\frac{c}{2}){\rm
ch}((\triangle(V))^b)
 {\rm
ch}(\Theta_2(T_CM,V_C,\xi_C,a,b))\right\}^{(4k)},\eqno(3.26)$$
Define virtual complex vector bundle $\widetilde{b}_r(T_{\bf C}M,
V_{\bf C},a,b)$ on $M$, $0\leq r\leq [\frac{k}{2}],$ via the
equality
$$\Theta_2(T_{C}M,V_C,\xi_C,a,b)\equiv \sum
_{r=0}^{[\frac{k}{2}]}\widetilde{b}_r(8\delta_2)^{k-2r}\varepsilon_2^r~~~~{\rm
mod}~q^{\frac{[\frac{k}{2}]+1}{2}}\dot
K(M)[[q^{\frac{1}{2}}]].\eqno(3.27)$$
 \noindent Define degree $4k-4$ differential forms
 $\widetilde{\beta}_r(T_{\bf C}M, V_{\bf C},,\xi_C,a,b)$ on $M$, $0\leq r\leq
[\frac{k}{2}],$ via the equality
$$\left\{\frac{e^{\frac{1}{24}E_2(\tau)[p_1(TM)-(a+2b)p_1(V)]}-1}{p_1(TM)-(a+2b)p_1(V)}\widehat{A}(TM)
{\rm ch}((\triangle(V))^b){\rm cosh}(\frac{c}{2}){\rm
ch}(\Theta_2(T_{C}M,V_C,\xi_C,a,b))\right\}^{(4k-4)}$$ $$\equiv \sum
_{r=0}^{[\frac{k}{2}]}\beta_r(8\delta_2)^{k-2r}\varepsilon_2^r~~~~{\rm
mod}~q^{\frac{[\frac{k}{2}]+1}{2}}\dot
\Omega^{(4k-4)}(M)[[q^{\frac{1}{2}}]].\eqno(3.28)$$
 Then similar to Theorem 3.1, we get\\

\noindent {\bf Theorem 3.4}
$$\left\{\frac{\widehat{A}(TM,\nabla^{TM}){\rm ch}((\triangle(V))^a)}{{\rm
cosh}^2(\frac{c}{2})}\right\}^{(4k)}$$
$$-\sum_{r=0}^{[\frac{k}{2}]}2^{(a-b)l+k-6r}\left\{\widehat{A}(TM,\nabla^{TM}){\rm cosh}(\frac{c}{2})
{\rm ch}((\triangle(V))^b){\rm ch}(\widetilde{b}_r)\right\}^{(4k)}$$
$$=[p_1(TM)-(a+2b)p_1(V)]\widetilde{{\cal
B}}(\nabla^{TM},\nabla^V,a,b),\eqno(3.29)$$
 {\it where}
$$\widetilde{{\cal
B}}(\nabla^{TM},\nabla^V,\nabla^\xi,a,b)=\sum_{r=0}^{[\frac{k}{2}]}2^{(a-b)l+k-6r}\widetilde{\beta_r}$$
$$-\left\{\frac{e^{\frac{1}{24}[p_1(TM)-(a+2b)p_1(V)]}-1}{p_1(TM)-(a+2b)p_1(V)}
\frac{\widehat{A}(TM,\nabla^{TM}){\rm ch}((\triangle(V))^a)}{{\rm
cosh}^2(\frac{c}{2})}\right\}^{(4k-4)}.\eqno(3.30)$$\\

\section{The anomaly cancellation formulas involving two complex
line bundles}

\quad Let $\xi, \xi'$ be two rank two real Euclidean vector bundle
with Euclidean connections $\nabla^\xi,~\nabla^{\xi'}$. Set
   $$\Theta_1(T_{C}M,V_C,\xi_C,\xi_C')=
   \bigotimes _{n=1}^{\infty}S_{q^n}(\widetilde{T_CM})\otimes
\bigotimes
_{m=1}^{\infty}\wedge_{q^m}(\widetilde{V_C}-2\widetilde{\xi_C})$$
$$~~~~~~~~\otimes \bigotimes _{r=1}^{\infty}\wedge
_{q^{r-\frac{1}{2}}}(\widetilde{\xi_C'})\otimes\bigotimes
_{s=1}^{\infty}\wedge _{-q^{s-\frac{1}{2}}}(\widetilde{\xi_C'}),$$
$$\Theta_2(T_{C}M,V_C,\xi_C,\xi_C')=\bigotimes _{n=1}^{\infty}S_{q^n}(\widetilde{T_CM})\otimes
\bigotimes _{m=1}^{\infty}\wedge_{q^{m}}(\widetilde{\xi_C'})$$
$$~~~~~~~~\otimes \bigotimes _{r=1}^{\infty}\wedge
_{q^{r-\frac{1}{2}}}(\widetilde{\xi_C'})\otimes\bigotimes
_{s=1}^{\infty}\wedge
_{-q^{s-\frac{1}{2}}}(\widetilde{V_C}-2\widetilde{\xi_C}),\eqno(4.1)$$
Set
$$P_1(\tau)=\left\{
\frac{e^{\frac{1}{12}E_2(\tau)[p_1(\xi)-p_1(\xi')]}{\widehat{A}(TM,\nabla^{TM})}{\rm
ch}(\triangle(V))}{{\rm cosh}^2(\frac{c}{2})}
 {\rm
ch}(\Theta_1(T_CM,V_C,\xi_C,\xi_C'))\right\}^{(4k)},\eqno(4.2)$$
$$P_2(\tau)=\left\{{\widehat{A}(TM,\nabla^{TM})}{\rm cosh}(\frac{c'}{2})
 {\rm
ch}(\Theta_2(T_CM,V_C,\xi_C,\xi_C'))\right\}^{(4k)},\eqno(4.3)$$
$$P_3(\tau)=\left\{\frac{e^{\frac{1}{12}E_2(\tau)[p_1(\xi)-p_1(\xi')]}-1}{p_1(\xi)-p_1(\xi')}
\right.$$
$$\left.\cdot{\widehat{A}(TM,\nabla^{TM})}{\rm cosh}(\frac{c'}{2})
 {\rm
ch}(\Theta_2(T_CM,V_C,\xi_C,\xi'_C))\right\}^{(4k)},\eqno(4.4)$$

Define virtual complex vector bundle $\overline{b}_r(T_{\bf C}M,
V_{\bf C},\xi,\xi')$ on $M$, $0\leq r\leq [\frac{k}{2}],$ via the
equality
$$\Theta_2(T_{C}M,V_C,\xi,\xi')\equiv \sum
_{r=0}^{[\frac{k}{2}]}\overline{b}_r(8\delta_2)^{k-2r}\varepsilon_2^r~~~~{\rm
mod}~q^{\frac{[\frac{k}{2}]+1}{2}}\dot
K(M)[[q^{\frac{1}{2}}]].\eqno(4.5)$$ Then
$$\overline{b}_0=(-1)^k{\bf C},~~\overline{b}_1=-24(-1)^kk+2\widetilde{\xi_C}+\widetilde{\xi_C}'
-\widetilde{V_{\bf C}}.\eqno(4.6)$$
 \indent Define degree $4k-4$ differential forms
 $\overline{\beta}_r(T_{\bf C}M, V_{\bf C},\xi,\xi')$ on $M$, $0\leq r\leq
[\frac{k}{2}],$ via the equality
$$\left\{\frac{e^{\frac{1}{12}E_2(\tau)[p_1(\xi)-p_1(\xi')]}-1}{p_1(\xi)-p_1(\xi')}\widehat{A}(TM)
{\rm cosh}(\frac{c'}{2}) {\rm
ch}(\Theta_2(T_{C}M,V_C,\xi,\xi'))\right\}^{(4k-4)}$$ $$\equiv \sum
_{r=0}^{[\frac{k}{2}]}\overline{\beta}_r(8\delta_2)^{k-2r}\varepsilon_2^r~~~~{\rm
mod}~q^{\frac{[\frac{k}{2}]+1}{2}}\dot
\Omega^{(4k-4)}(M)[[q^{\frac{1}{2}}]].\eqno(4.7)$$ It is easy to
calculate that
$$\overline{\beta}_0=(-1)^k\left\{\frac{e^{\frac{1}{12}
[p_1(\xi)-p_1(\xi')]}-1}{p_1(\xi)-p_1(\xi')}\widehat{A}(TM) {\rm
cosh}(\frac{c'}{2}) \right\}^{(4k-4)},\eqno(4.8)$$
$$\overline{\beta}_1=(-1)^k\left\{\frac{e^{\frac{1}{12}
[p_1(\xi)-p_1(\xi')]}-1}{p_1(\xi)-p_1(\xi')}\widehat{A}(TM) {\rm
cosh}(\frac{c'}{2}){\rm
ch}(-24k+2\widetilde{\xi_C}+\widetilde{\xi_C}' -\widetilde{V_{\bf
C}}) \right\}^{(4k-4)},\eqno(4.9)$$ We have
$$P_1(\tau)=2^{l}\left\{e^{\frac{1}{12}E_2(\tau)
[p_1(\xi)-p_1(\xi')]}\left(\prod_{j=1}^{2k}\frac{x_j\theta'(0,\tau)}{\theta(x_j,\tau)}\right)
\left(\prod_{\nu=1}^{l}\frac{\theta_1(y_\nu,\tau)}{\theta_1(0,\tau)}
\right)\right.$$
$$\left.\cdot
\left(\frac{\theta^2_1(0,\tau)}{\theta^2_1(u,\tau)}\frac{\theta_3(u',\tau)}{\theta_3(0,\tau)}
\frac{\theta_2(u',\tau)}{\theta_2(0,\tau)}\right)\right\}^{(4k)}
.\eqno(4.10)$$ Similarly,
$$P_2(\tau)+[p_1(\xi)-p_1(\xi')]P_3(\tau)=
\left\{e^{\frac{1}{12}E_2(\tau)
[p_1(\xi)-p_1(\xi')]}\left(\prod_{j=1}^{2k}\frac{x_j\theta'(0,\tau)}{\theta(x_j,\tau)}\right)\right.$$
$$\left.\cdot
\left(\prod_{\nu=1}^{l}\frac{\theta_2(y_\nu,\tau)}{\theta_2(0,\tau)}
\right)
\left(\frac{\theta^2_2(0,\tau)}{\theta^2_2(u,\tau)}\frac{\theta_3(u',\tau)}{\theta_3(0,\tau)}
\frac{\theta_1(u',\tau)}{\theta_1(0,\tau)}\right)\right\}^{(4k)}
.\eqno(4.11)$$ We assume that $p_1(TM)=p_1(V)$, then we have
$P_1(\tau)$ is a modular form of weight $2k$ over $\Gamma_0(2)$,
while $P_2(\tau)+[p_1(\xi)-p_1(\xi')]P_3(\tau)$ is a modular form of
weight $2k$ over $\Gamma^0(2)$ . Moreover, the following identity
holds,
$$P_1(-\frac{1}{\tau})=2^{l}\tau^{2k}(P_2(\tau)+[p_1(\xi)-p_1(\xi')]P_3(\tau)).\eqno(4.12)$$
So similar to the discussions in Section 3, we get\\

\noindent {\bf Theorem 4.1} {\it If $p_1(TM)=p_1(V)$, then}
$$\left\{\frac{\widehat{A}(TM,\nabla^{TM}){\rm ch}(\triangle(V))}{{\rm
cosh}^2(\frac{c}{2})}\right\}^{(4k)}$$
$$-\sum_{r=0}^{[\frac{k}{2}]}2^{l+k-6r}\left\{\widehat{A}(TM,\nabla^{TM}){\rm cosh}(\frac{c'}{2})
{\rm ch}(\overline{b}_r)\right\}^{(4k)}$$
$$=[p_1(\xi)-p_1(\xi')]\overline{{\cal
B}}(\nabla^{TM},\nabla^V,\nabla^\xi,\nabla^{\xi'}),\eqno(4.13)$$
 {\it where}
$$\widetilde{{\cal
B}}(\nabla^{TM},\nabla^V,\nabla^\xi,\nabla^{\xi'})=\sum_{r=0}^{[\frac{k}{2}]}2^{l+k-6r}\overline{\beta_r}$$
$$-\left\{\frac{e^{\frac{1}{12}[p_1(\xi)-p_1(\xi')]}-1}{p_1(\xi)-p_1(\xi')}
\frac{\widehat{A}(TM,\nabla^{TM}){\rm ch}(\triangle(V))}{{\rm
cosh}^2(\frac{c}{2})}\right\}^{(4k-4)}.\eqno(4.14)$$\\
\\
\noindent {\bf Corollary 4.2} {\it When ${\rm dim}M=4$, the
following identity holds}
$$\left\{\frac{\widehat{A}(TM,\nabla^{TM}){\rm ch}(\triangle(V))}{{\rm
cosh}^2(\frac{c}{2})}\right\}^{(4)}+2^{l+1}\left\{\widehat{A}(TM,\nabla^{TM}){\rm
cosh}(\frac{c'}{2})\right\}^{(4)}$$
$$=-2^{l-2}[p_1(\xi)-p_1(\xi')].\eqno(4.15)$$\\

\noindent {\bf Corollary 4.3} {\it When ${\rm dim}M=8$, the
following identity holds}
$$\left\{\frac{\widehat{A}(TM,\nabla^{TM}){\rm ch}(\triangle(V))}{{\rm
cosh}^2(\frac{c}{2})}\right\}^{(8)}-2^{l}\left\{\widehat{A}(TM,\nabla^{TM}){\rm
cosh}(\frac{c'}{2})\right\}^{(8)}$$
$$-2^{l-4}\left\{\widehat{A}(TM,\nabla^{TM}){\rm
cosh}(\frac{c'}{2}){\rm
ch}(2\widetilde{\xi_C}+\widetilde{\xi_C'}-\widetilde{V_C})
\right\}^{(8)}$$
$$
=[p_1(\xi)-p_1(\xi')]
\left\{\frac{e^{\frac{1}{12}[p_1(\xi)-p_1(\xi')]}-1}{p_1(\xi)-p_1(xi')}\right.$$
$$\left.\cdot
\left[2^{l}\widehat{A}(TM,\nabla^{TM}){\rm cosh}(\frac{c'}{2})-
\widehat{A}(TM,\nabla^{TM}){\rm cosh}(\frac{c'}{2}){\rm
ch}(2\widetilde{\xi_C}+\widetilde{\xi_C}' -\widetilde{V_{\bf
C}})\right.\right.$$
$$\left.\left.-\frac{\widehat{A}(TM,\nabla^{TM}){\rm ch}(\triangle(V))}{{\rm
cosh}^2(\frac{c}{2})}\right]\right\}^{(8)}.\eqno(4.16)$$

\section{The odd dimensional case}
 \quad In this section, Let $M$ be a (4k-1)-dimensional manifold. In the
 definition (3.22), we set $a=1,~b=0$ and $\Phi_1(\nabla^{TM},\nabla^V,\nabla^\xi,\tau)=\widetilde{Q_1}(\tau),$ and
$\Phi_2(\nabla^{TM},\nabla^V,\nabla^\xi,\tau)
=\widetilde{Q_2}(\tau)+[p_1(TM)-p_1(V)]\widetilde{Q_3(\tau)}.$
Applying the Chern-Weil theory, we can express $\Phi_1,~\Phi_2$ as
follows:

$$\Phi_1(\nabla^{TM},\nabla^V,\nabla^\xi,\tau)=2^le^{\frac{1}{24}E_2(\tau)[p_1(TM)-p_1(V)]}
{\rm
det}^{\frac{1}{2}}\left(\frac{R^{TM}}{4\pi^2}\frac{\theta'(0,\tau)}
{\theta(\frac{R^{TM}}{4\pi^2},\tau)}\right)$$ $$\cdot {\rm
det}^{\frac{1}{2}}\left(\frac{\theta_1(\frac{R^{V}}{4\pi^2},\tau)}{\theta_1(0,\tau)}\right)
{\rm det}^{\frac{1}{2}}\left(\frac{\theta_1^2(0,\tau)}
{\theta_1^2(\frac{R^{\xi}}{4\pi^2},\tau)}\frac{\theta_3(\frac{R^{\xi}}{4\pi^2},\tau)}{\theta_3(0,\tau)}
\frac{\theta_2(\frac{R^{\xi}}{4\pi^2},\tau)}{\theta_2(0,\tau)}\right);\eqno(5.1)$$
$$\Phi_2(\nabla^{TM},\nabla^V,\nabla^\xi,\tau)=e^{\frac{1}{24}E_2(\tau)[p_1(TM)-p_1(V)]}
{\rm
det}^{\frac{1}{2}}\left(\frac{R^{TM}}{4\pi^2}\frac{\theta'(0,\tau)}
{\theta(\frac{R^{TM}}{4\pi^2},\tau)}\right)$$ $$\cdot {\rm
det}^{\frac{1}{2}}\left(\frac{\theta_2(\frac{R^{V}}{4\pi^2},\tau)}{\theta_2(0,\tau)}\right)
{\rm det}^{\frac{1}{2}}\left(\frac{\theta_2^2(0,\tau)}
{\theta_2^2(\frac{R^{\xi}}{4\pi^2},\tau)}\frac{\theta_3(\frac{R^{\xi}}{4\pi^2},\tau)}{\theta_3(0,\tau)}
\frac{\theta_1(\frac{R^{\xi}}{4\pi^2},\tau)}{\theta_1(0,\tau)}\right);\eqno(5.2)$$

 Next we consider the transgression of
$\Phi_1(\nabla^{TM},\nabla^V,\nabla^\xi,\tau),~\Phi_1(\nabla^{TM},\nabla^V,\nabla^\xi,\tau).~$
about $\nabla^{\xi}$. Let $\nabla_1^{\xi},~\nabla_0^{\xi}$ be two
Euclidean connections on $\xi$ and
$B=\nabla_1^{\xi}-\nabla_0^{\xi},~\nabla_t^\xi=t\nabla_1^{\xi}+(1-t)\nabla_0^{\xi}$.
We
have\\
 $\Phi_1(\nabla^{TM},\nabla^V,\nabla_1^{\xi},\tau)
-\Phi_1(\nabla^{TM},\nabla^V,\nabla_0^{\xi},\tau)$ $$=
\frac{1}{8\pi^2}d\int_0^1
\Phi_1(\nabla^{TM},\nabla^V,\nabla_t^{\xi},\tau){\rm tr}\left[B
\left( \frac{\theta'_2(\frac{R_t^{\xi}}{4\pi^2},\tau)}
{\theta_2(\frac{R_t^{\xi}}{4\pi^2},\tau)}+\frac{\theta_3'(\frac{R_t^{\xi}}{4\pi^2},\tau)}
{\theta_3(\frac{R_t^{\xi}}{4\pi^2},\tau)}-2\frac{\theta_1'(\frac{R_t^{\xi}}{4\pi^2},\tau)}
{\theta_1(\frac{R_t^{\xi}}{4\pi^2},\tau)}\right)\right]dt.\eqno(5.3)$$
We define\\
$CS\Phi_1(\nabla^{TM},\nabla^V,\nabla_0^{\xi},\nabla_1^{\xi},\tau)$
$$:= \frac{\sqrt{2}}{8\pi^2}\int_0^1
\Phi_1(\nabla^{TM},\nabla^V,\nabla_t^{\xi},\tau){\rm tr}\left[B
\left( \frac{\theta'_2(\frac{R_t^{\xi}}{4\pi^2},\tau)}
{\theta_2(\frac{R_t^{\xi}}{4\pi^2},\tau)}+\frac{\theta_3'(\frac{R_t^{\xi}}{4\pi^2},\tau)}
{\theta_3(\frac{R_t^{\xi}}{4\pi^2},\tau)}-2\frac{\theta_1'(\frac{R_t^{\xi}}{4\pi^2},\tau)}
{\theta_1(\frac{R_t^{\xi}}{4\pi^2},\tau)}\right)\right]dt.\eqno(5.4)$$
which is in $\Omega^{\rm odd}(M,{\bf C})[[q^{\frac{1}{2}}]].$ Since
$M$ is $4k-1$ dimensional,
$\{CS\Phi_1(\nabla^{TM},\nabla^V,\nabla_0^{\xi},\nabla_1^{\xi},\tau)\}^{(4k-1)}$
represents an element in $H^{4k-1}(M,{\bf C})[[q^{\frac{1}{2}}]]$.
Similarly, we can compute the transgressed forms for
$\Phi_2,$ and define\\
$CS\Phi_2(\nabla^{TM},\nabla^V,\nabla_0^{\xi},\nabla_1^{\xi},\tau)$
$$:= \frac{1}{8\pi^2}\int_0^1
\Phi_2(\nabla^{TM},\nabla^V,\nabla_t^{\xi},\tau){\rm tr}\left[B
\left( \frac{\theta'_3(\frac{R_t^{\xi}}{4\pi^2},\tau)}
{\theta_3(\frac{R_t^{\xi}}{4\pi^2},\tau)}+\frac{\theta_1'(\frac{R_t^{\xi}}{4\pi^2},\tau)}
{\theta_1(\frac{R_t^{\xi}}{4\pi^2},\tau)}-2\frac{\theta_2'(\frac{R_t^{\xi}}{4\pi^2},\tau)}
{\theta_2(\frac{R_t^{\xi}}{4\pi^2},\tau)}\right)\right]dt,\eqno(5.5)$$
which also lies in $\Omega^{\rm odd}(M,{\bf C})[[q^{\frac{1}{2}}]]$
and its top component represents elements in $H^{4k-1}(M,{\bf
C})[[q^{\frac{1}{2}}]]$. Similar to Theorem 3.4 in [W], we have\\

 \noindent {\bf Theorem 5.1} {\it Let $M$ be a $4k-1$ dimensional manifold and
$\nabla^{TM}$ be a connection on $TM$ and $\xi$ be a two dimensional
oriented Euclidean real vector bundle with two Euclidean connections
$\nabla_1^{\xi}$,~$\nabla_0^{\xi}$, then we have
$\{CS\Phi_1(\nabla^{TM},\nabla^V,\nabla_0^{\xi},\nabla_1^{\xi},\tau)\}^{(4k-1)}$
is a modular form of weight $2k$ over $\Gamma_0(2)$;
$\{CS\Phi_2(\nabla^{TM},\nabla^V,\nabla_0^{\xi},\nabla_1^{\xi},\tau)\}^{(4k-1)}$
is a modular form of weight $2k$ over $\Gamma^0(2);$
 The following equalities hold,}
$$\{CS\Phi_1(\nabla^{TM},\nabla^V,\nabla_0^{\xi},\nabla_1^{\xi},\tau)\}^{(4k-1)}
=(2\tau)^{2k}\{CS\Phi_2(\nabla^{TM},\nabla^V,\nabla_0^{\xi},\nabla_1^{\xi},\tau)\}^{(4k-1)},$$

\indent Thus by Theorem 5.1, we can get a cancellation formula
similar to Corollary 3.5 in [W].\\

 \noindent{\bf Acknowledgement}~The
work of the second author was supported by NSFC No.10801027 and Fok
Ying Tong Education Foundation No. 121003.
\\

\noindent {\bf References}\\

\noindent [AW] L. Alvarez-Gaum\'{e}, E. Witten, Graviational
anomalies, {\it Nucl.Phys.} B234 (1983), 269-330.\\
 \noindent [C] K. Chandrasekharan, {\it Elliptic
Functions},
Spinger-Verlag, 1985. \\
\noindent [CH1] Q. Chen, F. Han, Modular invariance and twisted
anomaly cancellations of characteristic numbers, {\it Trans. Amer.
Math. Soc.} 361 (2009), 1463-1493 \\
\noindent [CH2] Q. Chen, F. Han, Elliptic genera, transgression and
loop space Chern-Simons form, {\it Comm. Anal. Geom.}, 17
(2009),:73-106.\\
\noindent [HLZ] F. Han, K. Liu, W. Zhang, Modular Forms and
Generalized Anomaly Cancellation Formulas, arXiv:1109.2494.\\
 \noindent [HZ1]F. Han, W. Zhang, ${\rm
Spin}^c$-manifold and elliptic genera,
{\it C. R. Acad. Sci. Paris Serie I.,} 336 (2003), 1011-1014.\\
\noindent [HZ2] F. Han, W. Zhang, Modular invariance, characteristic
numbers and eta Invariants,  {\it J.
Diff. Geom.} 67 (2004), 257-288.\\
\noindent [L] K. Liu, Modular invariance and characteristic
numbers. {\it Commu. Math. Phys.} 174 (1995), 29-42.\\
\noindent [W] Y. Wang, Transgression and twisted anomaly
cancellation formulas on odd dimensional manifolds, {\it J. Geom.
Phys.},60 (2010), 611-622.\\
 \noindent [Z] W. Zhang, {\it Lectures on Chern-weil Theory
and Witten Deformations.} Nankai Tracks in Mathematics Vol. 4, World
Scientific, Singapore, 2001.\\

\indent{ Center of Mathematical Sciences, Zhejiang University
Hangzhou Zhejiang 310027, China and Department of Mathematics,
University of California at Los Angeles, Los Angeles CA 90095-1555,
USA\\} \indent  Email: {\it liu@ucla.edu.cn; liu@cms.zju.edu.cn}\\

 \indent{  School of Mathematics and Statistics,
Northeast Normal University, Changchun Jilin, 130024 China }\\
\indent E-mail: {\it wangy581@nenu.edu.cn; }\\

\end {document}